\documentclass[12pt]{amsproc}
\usepackage{amsfonts,amssymb,amscd,amsthm,amsmath}
\usepackage{graphpap,xypic}

 

\newtheorem*{theorem*}{Theorem}

\newtheorem*{lemma*}{Lemma}

\newtheorem*{proposition*}{Proposition}

\theoremstyle{definition}

\newtheorem*{definition*}{Definition}

\theoremstyle{remark}

\newtheorem*{remark*}{Remark}
\newcommand{\mt}[1]{\operatorname{#1}}

\newcommand{\QQ}{{\mathbb Q}}
\newcommand{\ZZ}{{\mathbb Z}}

\newcommand{\OO}{{\mathcal O}}

\newcommand{\PP}{{\mathbb P}}

\newcommand{\HH}{{\mathcal{H}}}

\newcommand{\FFF}{{\mathbb F}}

\newcommand{\Pic}{\mt{Pic}}

\newcommand{\mult}{\mt{mult}}

\title{On the Noether--Fano inequalities}
\author{\Large V.~A.~Iskovskikh}

\begin{document}
\maketitle

\parshape=1
3cm 10cm \noindent {\small \quad \quad \quad
\quad\quad\quad\quad\quad\quad\quad {\bf Abstract}\newline I give
a survey of Noether-Fano inequalities in birational geometry,
starting with the original Noether inequality and up to the
modern approach of Log Minimal Model program. The paper is based
on my talk at the Fano conference in Torino in October 2002. This
is the revised version: an erroneous reference in my paper
published in the Proceedings of Fano conference is corrected.}
\vspace{1cm}

\section{The Noether inequality}
\subsection{} The general theory of birational maps of plane has been developed
by Cremona since 1863. The factorization conjecture of such maps
into quadratic ones was given by Klebsh in 1869. He proved it for
the maps of degree $n\le 8$. M.~Noether also announced the proof
of this conjecture in 1869 and showed that it follows from the
inequality
\begin{equation}\label{eq1}
\nu_1+\nu_2+\nu_3>n
\end{equation}
for three maximal multiplicities of base points of linear system
giving the birational map. He proved this inequality in 1871.
\par
Independently Rosanes discovered this inequality and gave his own
proof by induction.
\\

\par
However Noether's factorization theorem does follow from this
inequality very non-trivially. If there exists a quadratic map
with base points satisfying (\ref{eq1}) then its composition with
the initial birational map has a lesser degree, but the existence
of a required quadratic untwisting map is a hard problem. A lot of
papers of different authors were devoted to it. It is accepted
that the first complete proof of Noether's theorem was obtained by
Castelnuovo in 1901 and the most clear proof was obtained by
Alexander in 1916 (see \cite{Hu}).
\\

\par
The classical proofs of Noether's inequality (\ref{eq1})
immediately follow from two basic equalities

\begin{gather}
\sum \nu_i^2=n^2-1\\
\sum \nu_i=3n-3,
\end{gather}
where $n$ is the degree of linear system giving Cremona
transformation $\chi\colon\PP^2\dashrightarrow\PP^2$ and $\nu_i$
are the multiplicities of base points including infinitely close
points of this linear system.

\subsection{} The birational and geometrical sense of Noether's inequality
was understood later and it is close to the question of birational
invariance of adjunction termination. By modern proof of Noether's
theorem according to Castelnuovo (see [I--R] preface to Hudson's
book and 1.5 below) it is not necessarily to find the triple of
base points with maximal multiplicities satisfying (\ref{eq1}). It
is enough to find one maximal singularity for the linear systems
on $\PP^2$ and $\FFF_N$ only, where
$\FFF_N=\PP_{\PP^1}(\OO+\OO(N))$ is a standart linear surface.

\subsection{Noether lemma}\label{ref1}
\emph{Let us consider the main ideas of the classical proof. Let
$F$ be a surface $\PP^2$ or $\FFF_N$ and $\chi\colon
F\dashrightarrow F'=\PP^2$ be a birational map. Let $\HH'$ be a
linear system of lines on $\PP^2$ and $\HH=\chi^{-1}_*(\HH')$ be a
proper transform on $F$. Write $H\sim \alpha s+\beta f$, where
$H\in \HH$ is a general divisor, $\alpha,\beta\in \ZZ$, $\alpha\ge
1$, $f$ is a fiber of ruled surface, $s$ is an exceptional section
of $F=\FFF_N$ or $s$ is a line on $F=\PP^2$. It is better to write
it in another basis $\{-K_F,f\}$ of the space $\Pic(F)\otimes\QQ$
: $H'\sim_{\QQ}-aK_F+bf$, where $a=\alpha/3$, $b=0$ for $F=\PP^2$,
and $a=\alpha/2$, $b=\beta-\frac{N+2}2\alpha$ for $F=\FFF_N$.}
\par
\emph{Then, if $\chi$ is not isomorphism then}
\begin{enumerate}
\item $\HH$ \emph{has a base point $P$ with a multiplicity $\nu>a$; or
\item $F=\FFF_N$ and $b<0$.}
\end{enumerate}

\subsection{Remark} In the notations of Log-MMP both these
statements are equivalent to the following:
\begin{description}
 \item if $K_F+\frac1aH$ is
canonical and $b\ge0$ then $\chi$ is an isomorphism.
\end{description}

\begin{proof} Consider a resolution
\\
\begin{center}
\begin{picture}(80,46)(0,0)
\large{ \put(14,10){$F$} \put(40,40){\vector(-1,-1){18}}
\put(42,42){$Z$} \put(50,40){\vector(1,-1){18}}
\put(64,10){$F'=\PP^2,$} \put(24,34){$\sigma$}
\put(62,34){$\varphi$} \put(38,10){$\dashrightarrow$}}
\end{picture}
\end{center}
$\HH_Z=\varphi^*\HH'=\sigma^{-1}_*\HH$,
$K_Z+\frac1aH_Z=\sigma^*(K_F+\frac1aH)+\sum_i(1-\frac{\nu_i}a)\sigma^*_iE_i$,
where $\sigma=\sigma_r\circ\ldots\circ\sigma_1$ and $E_i$ are
exceptional, $H_Z\in\HH_Z$, $H\in\HH$, $H'\in\HH'$ are general
divisors. Applying $\varphi_*$ we get

\begin{gather*}
(-3+(1/a))H'\sim
K_{\PP^2}+\frac1aH'\sim\varphi_*\sigma^*(bf)+\\
+ \sum(1-\frac{\nu_i}a) \varphi_*\sigma_i^*E_i.
\end{gather*}
If $\nu_i\le a$ for all $i$ and $b\ge 0$ then $(-3+\frac1a)H'$ is
an effective divisor. It can happen only if $a=\frac13$, i.e.
$\alpha=1$, $F=\PP^2$ and $\HH=\HH'$ and hence $\chi$ is an
isomorphism.
\end{proof}

\subsection{Noether-Castelnuovo theorem}
\begin{theorem*} Any birational map $\PP^2\dashrightarrow\PP^2$ is the composition
of the following elementary transformations (links):
\begin{enumerate}
\item[A).] a blow-up of a point
$\sigma^{-1}\colon\PP^2\dashrightarrow \FFF_1;$ \item[B).] an
elementary transformation $\varepsilon\colon\FFF_N\dashrightarrow
\FFF_{N\pm 1};$ \item[$A^{-1}$).] a contraction of exceptional
curve $\sigma\colon\FFF_1\to\PP^2$; \item[C).] a biregular
involution $\tau\colon\FFF_0\to \FFF_0$.
\end{enumerate}
\end{theorem*}
\begin{proof}
In the previous notations, if
$\chi\colon\PP^2\dashrightarrow\PP^2$ is not isomorphism then by
\ref{ref1} a linear system $\HH$ has a base point (perhaps,
infinitely close one) with a multiplicity $\nu>a$. Since the
multiplicities don't increase on the resolution then such a point
exists on $\PP^2$. Apply link A)
$\sigma^{-1}\colon\FFF_1\dashrightarrow\PP^2$ with the senter in
this point. Let $\HH_1=\sigma^{-1}_*\HH$. Then
\begin{gather*}
H_1\sim-\frac{(3a-\nu)}2K_{\FFF_1}+\frac{3(\nu-a)}2f_1=\\
=-a_1K_{\FFF_1}+b_1f_1.
\end{gather*}
where $H_1\subset\HH_1$ is a general divisor. Applying the links
of type B) we can untwist all maximal singularities of
multiplicity $>a_1$ and we obtain a linear system $\HH_N\sim
-a_1K_{\FFF_N}+b_Nf$ on the surface $\FFF_N$ for some $N$. By
\ref{ref1} $b_N<0$. Then there are two cases $N=0$ or $N=1$ only.
Indeed, $H_N$ doesn't have the fixed components and hence it is a
nef and big divisor. By index theorem, if $H_N\cdot C=0$ for some
curve $C$ then $C^2<0$, i.e. $C=s_N$. Also $0\le s_N\cdot
H_N<-a_1K_{\FFF_N}\cdot s_N=2-N$, i.e. $N\le 1$. A divisor
$-K_{\FFF_N}$ is ample if and only if $N=0,1$.

If $N=1$ let us apply link $A^{-1}$) $\sigma\colon\FFF\to\PP^2$.
Then $\HH'=\sigma_*\HH_1\sim a'(-K_{\PP^2})$ with $a'=a_1+\frac13
b_{N}=a+\frac{a-\nu}2+\frac13b_{N}<a$, since $\nu>a$ and
$b_{N}<0$.

In the case $N=0$ we have $H_{N}\sim
-a_1K+b_Nf_N=2a_1S_0+(2a_1+b_0)f_0=-(2a_1+b_0)K_{F_0}-b_0s_0$, so
the statement (1) of (1.3) is true if it is considered on another
structure of $\PP^1$-fibration $\FFF_0\to \PP^1$, since
$\FFF_0\ncong\PP^2$ and $-b_0>0$. So after links $\tau$ and $B)$
we fall into the previous situation , but with a smaller
coefficient of $-K_F$ (see also \cite{I-R} and \cite{I2}).
\end{proof}

\subsection{Remark}
In fact, just these very reasonings and not very clear notes of
Sarkisov have inspired M. Reid [R] to formulate the program about
decomposition into elementary links of four types of birational
maps between Mori-fibred spaces $\{\phi:X\to S\}$. He named it
Sarkisov program. In dimension 2 the existing of decomposition and
finiteness of this algorithm were just shown in (1.4). However it
could be seen as a particular case of the general theorem of
Corti, which says the existence and finiteness of the algorithm in
Sarkisov-Reid program in dimension at most three (see the
following paragraph).

\section {Fano inequality and Sarkisov-Reid program}
\subsection{}
Studying of birational transformation of $\PP^3$ was started in
Cayley (1869-1870), Noether (1870-1871) and Cremona (1871-1872)
papers (see, for example, [Hu]). However only some particular
examples were considered and no general theory (as in dimension 2)
was created. In the early of twenty century there was an
understanding (through the birational invariant of the termination
of adjunction) that for Cremona transformation of $\PP^3$ there is
an analog of the Noether inequality: there is  a maximal
singularity: either a basis curve of degree more than $\frac n4$,
or a point with a multiplicity more than $\frac n2$, where $n$ -
is degree of the linear system $\HH$ in $\PP^3$ which define this
transformation and that the joint degree of maximal curves are
$\le 15$.

I do not know whether the existing of infinitely close maximal
curves for three-dimensional Cremona transformations was analyzed,
but Fano in the paper ([Fa], 1915) essentially used  the analog of
the Noether inequality , including the infinitely closed maximal
curve (only on the first blow up of a point) for his studying of
birational characteristics of three-dimensional quartic
$V_4\subset \PP^4$ and the full intersection $V_{2,3}\subset
\PP^5$.

\subsection{
Classical Fano inequality, 1915}
\begin{lemma*}
Let $V$ and $V'$ be  smooth Fano threefolds, $\rho(V)=\rho(V')=1$
where $\rho$ is Picard number. Let $\chi:V\dashrightarrow V'$ be a
birational map, $H'$ be the positive generator in $\Pic(V')$, $H$
be the positive generator in $\Pic(V)$. Let
$\mathcal{M}=\chi^{-1}_*(|H'|)\subset|nH|$ be proper transform of
the linear system $\HH'=|H'|$. Then, if $\chi$ is not an
isomorphism, then $\mathcal{M}$ has a maximal singularity of one
of the following types:
\begin{enumerate}
\item curve $C\subset V$, $\mult_C \mathcal{M}>\frac nr, \ \deg
C=CH<r^2H^3$, where $r$ ia an index of $V$, i.e. $-K_V\sim rH$;
\item point $P\in V, \ \mult_P \mathcal{M}>\frac {2n}r$;

\item infinitely close curve $B^*\subset
V^*\stackrel{\sigma}{\longrightarrow}V$, $\mult_B^*
\mathcal{M}_{V^*}>\frac nr$, $\mult_{\sigma(B^*)}\mathcal{M}
>\frac{n}r$, where $\sigma(B^*)$ is a point on $V$.
\end{enumerate}
\end{lemma*}

This lemma can be proved in the similar way as  lemma 1.3 did.
Generalizations on smooth Fano varieties and conic bundles are
made in [I-M], [Isk], [I-P], where these inequalities were named
Noether-Fano inequalities.

The first generalization onto singular varieties was formulated by
M. Reid within the framework of the Sarkisov-Reid program for
Mori-fibred spaces with $\QQ$-factorial varieties with terminal
singularities and was proved by Corti [Co] for treefolds (see also
[Ma]).

\subsection{Mori fibration}
\begin{definition*}
Let $X$ be a projective $\QQ$-factorial variety with terminal
singularities. The morphism $\varphi:X\to S$ is called Mori
fibration if
\begin{enumerate}
\item $\dim S<\dim X$, $\varphi_*\OO_X=\OO_S$;
\item $\rho(X/S)=1$;
\item $-K_X$ is $\varphi$-ample.
\end{enumerate}
\end{definition*}

Properties 1)-3) mean that $\varphi:X\to S$ is an extremal
contraction of a fibred type, $S$ is normal $\QQ$-factorial
variety.

\subsection{Sarkisov-Reid program}

This is a program (which is still hypothetical in dimension more
then three) of decomposition of every birational map between Mori
fibred spaces
\[
 \mbox{\definemorphism{birto}\dashed \tip \notip
\diagram
X \dto_\varphi \rrbirto^\chi& &\,X' \dto^{\varphi'} \\
 S & & S' \enddiagram}
 \]
(which is not necessarily preserve fibration structures) into the
finite composition of elementary links -- commutative diagrams of
one of the following four types

 \[
 \begin{array}{lclc}
 \raise-1.4cm\hbox{Type I:} &
\mbox{
\definemorphism{birto}\dashed \tip \notip \diagram
 &Z\dlto_\sigma\rbirto^{\psi}&X_2 \ar[dd] \\
X_1\dto_{\varphi_1} & & \dto\\
S_1  & &\ar[ll]^{\alpha} S_2\enddiagram}
\end{array}
 \]

where $\varphi_1,\ \varphi_2$ are Mori fibrations, $\sigma$ is a
divisorial extremal contraction, $\psi$ is a sequence of Mori
flips, anti-flips and flops, $\alpha$ is a surjective map with
connected fibres. (This is actually a high dimensional analog of
the link of the type A) in 1.4);

 \[
 \begin{array}{lclc}
 \raise-1.4cm\hbox{Type II:} &
\mbox{\definemorphism{birto}\dashed \tip \notip \diagram
 &Z_1\dlto_{\sigma_1}\rbirto^\psi&Z_2 \drto^{\sigma_2}& \\
X_1\drto_{\varphi_1}& & &X_2\dlto^{\varphi_2}\\
 &S_1\rdouble &S_2 & \enddiagram}
 \end{array}
 \]

 where $\varphi_1,\ \varphi_2$ are extremal divisorial
contractions Mori,  $\psi$ is a sequence of Mori flips, antiflips
and flops. (This is actually a high dimensional analog of the link
of the type B) in 1.4);

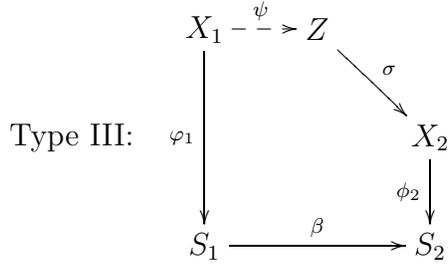
\[
 \begin{array}{lclc}
 \raise-1.4cm\hbox{Type III:} &
\mbox{\definemorphism{birto}\dashed \tip \notip \diagram
 X_1\ar[dd]_{\varphi_1}\rbirto^{\psi}&Z \ar[dr]^\sigma& \\
 & & X_2\dto_{\phi_2}\\
S_1 \ar[rr]^\beta & & S_2\enddiagram}
\end{array}
 \]
 This link is opposite to the one of type (I) (a higher dimensional
 analog of link of the type $A^{-1}$);
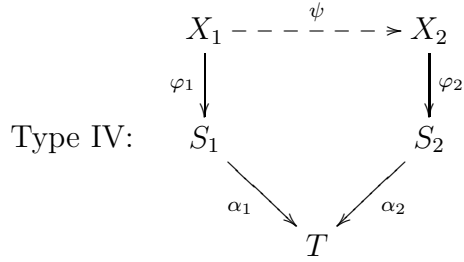
\[
 \begin{array}{lclc}
\raise-1.4cm\hbox{Type IV:} & \mbox{\definemorphism{birto}\dashed
\tip \notip \diagram
X_1\rrbirto^{\psi}\dto_{\varphi_1}& &X_2 \dto^{\varphi_2} \\
S_1 \drto_{\alpha_1} & &S_2 \dlto^{\alpha_2}\\
 &T& \enddiagram}
 \end{array}
 \]

 where $\varphi$ is a sequence of logflips, $\alpha_1,\
\alpha_2$ -- surjective morphisms with connected fibres (this is
actually a higher dimensional analog of link of the type C).

\subsection{}
To construct the algorithm of decomposition we need a numerically
ordered characteristic, which is decreasing during the
decomposition. We can use $\deg (\chi,\HH')$ as such a
characteristic. It is define by the map $\chi$ and previously
fixed on the whole process of decomposition very ample linear
system $\HH'=|-\mu' K_{V'}+\varphi'^*A'|$ in the following diagram
 \[\large{
 \xymatrix{
 (2.1)&\hbox{$\HH\sim-\mu K_X+\varphi^*A$,} &
X\dto_{\varphi}\ar@{-->}[r]^{\chi}&X' \dto^{\varphi'}&
\hbox{$\HH'=|-\mu'K_{X'}+\varphi'^*A'|$}
 \\
&\hbox{$\mu\in\QQ_{>0},A\in\Pic S$}&S&S'&\hbox{$A'\in\Pic S',\
\mu'\in\ZZ$}& & }
 }\]
where $\HH=\HH_X$ -- proper transform of the linear system $\HH'$.

By definition we have $\deg(\chi,\HH')=(\mu,\lambda,e),\
\mu\in\QQ_{>0},\ \lambda\in\QQ_{\ge0},\ e\in\ZZ_{\ge0}$, is
lexicographically ordered triple, where $\mu$ -- as in the diagram
2.1.

$\lambda=\frac1c$, where $c:=\max\{t\in\QQ_{>0}|K_X+tH_X\  is\
canonical\},\ H_X\in \HH$ -- is a general divisor, in other words,
$c$ is canonical threshold for log-pair $(X,H_X)$.

If
\\
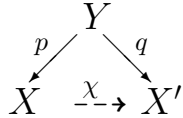
\begin{center}
\begin{picture}(80,46)(0,0)
\large{ \put(14,10){$X$} \put(40,40){\vector(-1,-1){18}}
\put(42,42){$Y$} \put(50,40){\vector(1,-1){18}} \put(64,10){$X'$}
\put(24,34){\footnotesize{$p$}} \put(62,34){\footnotesize{$q$}}
\put(38,10){$\dashrightarrow$} \put(42,18){\footnotesize{$\chi$}}
}
\end{picture}
\end{center}

is a common resolution for $\chi$, such that $p$ is log-resolution
for pair $(X,H_X)$ and
\begin{gather*}
K_Y=p^*K_X+\sum a_k E_k,\\
q^*H_{X'}=p^*H_{X'}-\sum b_kE_k,
\end{gather*}
where $E_k$ are all $p$-exceptional divisors, then
\begin{equation*}
\lambda=\max\left\{\frac{b_k}{a_k}\right\}
\end{equation*}
The property to be canonical for the pair $(X,cH_X)$ means that
inequalities for discrepancies $a(E,X,cH_X)\ge 0$ for all
exceptional divisors $E$ over $X$, i.e. $(a_k-cb_k)\ge0$ for every
$k$, and $c=\min\{\frac{a_k}{b_k}\}$, are true.

At the end, we have $e$ is number $\sharp\left\{E_k|\lambda
a_k-b_k=0\right\}$ of "maximal singularities".

\subsection{Remark}
In the case of smooth surfaces, as in 1.4, $\lambda$ is
multiplicity of maximal singularity of $\HH_X$. In the dimension
three and higher this connection with maximal multiplicity is not
quite straightforward. It is connected with the fact that for
smooth threefold discrepancies depend not only on incidences in
the graph of the exceptional divisors, but also on the weight 2 or
1 which appear with exceptional divisor when we blow up point or
curve.
\\

As to  singular (terminal) points, the situation is much more
complicated, because $a_k,\ b_k\in\QQ_{\ge0}$ are rational numbers
and decreasing of the degree $\deg (\chi, \HH')$, in particular
the process may be infinite. However, for the algorithm to be
finite one needs the break of descending sequens of
lexicographically ordered triples $(\mu,\lambda,e)$.

The similar way as we have in dimension 2, the main ingredient in
the process of decomposition is Noether-Fano inequality -- the
criterion of stopping of the decomposition.

\subsection{Noether-Fano criterion}
\emph{In the previous notations
\[
 \mbox{\definemorphism{birto}\dashed \tip \notip
\diagram
X \dto_\varphi \rrbirto^\chi& &\,X' \dto^{\varphi'} \\
 S & & S' \enddiagram}
 \]
$\chi$ is an isomorphism between two Mori-fibred spaces with
$\QQ$-factorial terminal singularities if $\lambda\le\mu$ (i.e the
linear system $\HH_X$ has no maximal singularities) and divisor
$K_X+\frac1\mu H_X$ is nef.}

The proof of the theorem is similar to the proof  of the classical
case (see [Co], [Ma]). We consider a common resolution
$X\stackrel{p}{\leftarrow}Y\stackrel{q}{\to}X'$ and study
intersections of divisors $(K_X+\frac 1{\mu'} H_X)$ and
$(K_{X'}+\frac 1\mu H_{X'})$ on $Y$ with generic curves in fibres
of morphisms $\varphi'p$ and $\varphi q$ using Negativity lemma.
In contrast to the previous reasoning with effective divisors
(geometrical case), here we use an intersection theory (nef case)
which is easier and more convenient to work within log-MMP. In
geometrical formulation the divisor $K_X+\frac 1\mu H_X$ (when
$\lambda\le\mu$) has to be effective (or quasi-effective).

\section{Generalizations}
\subsection{}
In the paper [Br-Ma] (see also [Ma]) log-variant of Sarkisov-Reid
program in the cathegory of $\QQ$-factorial Kawamata log-terminal
pairs $(X,B)$ with log-MMP relation is studied. \pagebreak
\subsection{Log-MMP relation}
\begin{definition*}
A finite number of projective log-pairs $(X_i,B_i)$, $i=1,\dots,k$
with only $\QQ$-factorial and klt singularities are said to be
log-MMP related iff there exists a log-pair $(Y,B_Y)$ with
nonsingular projective $Y$ and a boundary $\QQ$-divisor $B_Y$ with
only normal crossings such that all log-pairs $(X_i,B_i)$ are
obtained from $(Y,B_Y)$ via log-MMP.
\end{definition*}
For klt pairs of the fibred type $(X,B)\overset{\phi}{\to} S$ and
$(X',B_{X'})\overset{\phi'}{to}S'$ which are log-MMP related and
for a birational map between them
\[\large{
\xymatrix{
\hbox{(3.1)}&\hbox{$\HH\sim-\mu(K_X+B)+\varphi^*A$,}
 &(X,B)\dto_{\varphi}\ar@{-->}[r]^{\chi}&\hbox{$(X',B_{X'})$} \dto^{\varphi'} \\
& \hbox{$\HH'\sim-\mu'(K_{X'}+B_{X'})+\varphi'^*A'$},&S& S' }}
 \]
we can define a degree $(\mu, \lambda,e)$, where $\lambda$ is a
maximal multiplicity of extremal ray for $K_X+B_X+cH_X$ on a good
common log-resolution $(Y,B_Y)\to(X,B_X),\ (Y,B_Y)\to(X',B_{X'})$.

\subsection{Noether-Fano criterion for the log Sarkisov-Reid
program with klt singularities}
\begin{theorem*}
In the diagram(3.1) $\chi$ is an isomorphism if $\lambda\le\mu$
(i.e. $K_X+B_X+\frac 1\mu H_X$ is canonical) and $K_X+B_X+\frac
1\mu H_X$ is nef. In the geometrical variant the condition to be
nef is replaced to be effective.
\end{theorem*}
The most general variant of Noether-Fano inequality was proposed
by Shokurov and Cheltsov which is reduced to the statement about
the uniqueness of the canonical model for log-pair $(X,B)$ and
birational invariance of log-kodaira dimension.
\subsection{Canonical model}
\begin{definition*}
A pair $(V,B_V)$ is called canonical model for $(X,B)$ if there is
a birational map $\psi:X\dasharrow V$ such that
$(V,B_V)=(V,\psi(B))$ is canonical and divisor $K_V+B_V$ is ample.
\end{definition*}
\begin{proposition*}[Shokurov]
If a canonical model for pair $(X,B)$ exists then it is unique.
\end{proposition*}

\subsection{Iitaka map and Kodaira dimension}
\begin{definition*}
For pair $(X,B)$ we consider a birational map $\alpha:Y\dasharrow
X$ such that log pair $(Y,\phi^{-1}(B))=(Y,B_Y)$ is canonical. The
rational map
$\varphi=\varphi_n:(Y,B_Y)\dasharrow(Z,B_Z)=(\varphi_n(Y),\varphi_n(B_Y))$
define by the linear system $|n(K_Y+B_Y)|,\ n\gg 1$, is called
Iitaka map. By the Kodaira dimension we mean the number $\varkappa
(X,B):=\dim(Z,B_Z))$, if $|n(K_Y+B_Y)|\ne\emptyset$ for some $n$,
otherwise $\varkappa(X,B)=-\infty$.
\end{definition*}
\begin{proposition*}[Shokurov]
The map $\varphi$ and $\varkappa(X,B)$ not  depend on the
birational map $\alpha$.
\end{proposition*}

Such a generalisation allow us to study not only birational maps
between Mori-fibred spaces but birational maps between so-called
$K$-trivial bundles, in sense of [Ch1], i.e. those bundles the
general fibre of which has the Kodaira dimension 0.

This was showed in [Ch2] for  general smooth hypersurfaces of
degree $N$ in $\PP^N,\ N\ge4$.

\subsection{}
\begin{theorem*}[Ch2]
Let $X=X_N\subset\PP^N$ be a generic smooth hypersurfaces of
degree $N\ge4$. Then $X$ is not birational to fibrations, whose
generic fiber has Kodaira dimension 0, except for fibrations
induced by projections from $(N-2)$-dimensional linear subspace
in $\PP^N$.
\end{theorem*}
\begin{remark*}
In the first version of this paper [Isk3] the theorem on
birational geometry of generic Fano hypersurfaces
$X_N\subset{\mathbb P}^N$ was erroneously attributed to
I.Cheltsov (Sec. 3.6 of my paper).

The fact that a generic Fano hypersurface $X_N\subset{\mathbb
P}^N$ is birationally superrigid was proved by A.V.Pukhlikov and
published in [P1].

The claim, formulated in my paper as part (2) of the Theorem of
Sec. 3.6, is an immediate consequence of the superrigidity,
stated in [P1, Sec. 2].

The claim, formulated in my paper as part (1) of the Theorem of
Sec. 3.6, follows from the superrigidity in an elementary way.

Both statements 1) and 2) is proved for any smooth hypersurface
$X_N\subset \PP^N$ by Pukhlikov [Pu1] for $N\geq 6$.
\end{remark*}

The main idea of the proof of the theorem  is as follows. Let
$\HH$ be our linear system with no fixed components,
$\HH\subset|-nK_X|$, $CS(X,\frac1n H)$ be the locus of canonical
singularities for the pair $(X,\frac1n H)$, where $H\in\HH$ is a
general divisor. It was proved in [P1] that the pair $(X,\frac1n
H)$ is canonical. Using this fact, Cheltsov shows that either
$CS(X,\frac1n H)=\emptyset$ or $CS(X,\frac1n H)=(X\bigcap L)$ for
some $L\simeq\PP^{N-2}$.

We have if $CS=\emptyset$ (i.e. $(X,\frac 1n H)$ is terminal) then
$(X,(\frac 1n+\varepsilon)H)$ is still terminal for a small
$0<\varepsilon\ll 1$. Hence the divisor $K_X+(\frac
1n+\varepsilon)H\sim-\varepsilon K_X$ ia ample. The uniqueness of
the canonical model implies that $\HH$ defines an isomorphism.

If $CS\neq\emptyset$, then the Kodaira dimension
$\varkappa(X,H)=1$ and $\HH$ defines a birational map
$\varphi_{\HH}:X\dasharrow \PP^1$ which is birationally equivalent
to a projection bundle $\varphi_L:X\dasharrow \PP^1$ from a linear
space $L$. So $\varphi_{\HH}$ and $\varphi_L$ are birationally
equivalent $K$-trivial rational fibrations.

 \end{document}